\documentclass[12pt]{article}

\usepackage{fancyhdr}
\usepackage{amsmath}
\usepackage{amsthm}
\usepackage{amssymb}
\usepackage{lastpage}
\usepackage{hyperref}
\usepackage{graphicx}
\usepackage{datetime}
\usepackage{pdfpages}
\usepackage{caption}
\usepackage{booktabs}
\usepackage{float}
\usepackage{subcaption}
\usepackage{tikz}
\usepackage{natbib}
\usepackage{url}
\usepackage{xcolor}
\usepackage{authblk}
\usepackage{listings}

\usetikzlibrary{calc, arrows.meta}

\hypersetup{
    colorlinks,
    citecolor=blue
}

\newcommand{\abs}[1]{\left\vert#1\right\vert}
\newcommand{\set}[2]{\left\{\left. #1 \,\right| #2 \right\}}
\newcommand{\injr}[1]{\text{inj}(#1)}
\newcommand{\norm}[1]{\left\lVert#1\right\rVert}
\newcommand{\inpr}[2]{\left\langle #1, #2 \right\rangle}

\newtheorem{theorem}{\underline{Theorem}}[section]

\newtheorem{lemma}[theorem]{\underline{Lemma}}
\newtheorem{corollary}[theorem]{\underline{Corollary}}

\theoremstyle{definition}
\newtheorem{definition}[theorem]{\underline{Definition}}
\newtheorem{remark}[theorem]{\underline{Remark}}

\lstdefinestyle{mypython}{
    language=Python,
    basicstyle=\ttfamily\footnotesize,
    keywordstyle=\color{blue}\bfseries,
    commentstyle=\color{gray}\itshape,
    stringstyle=\color{red},
    numbers=left,
    numberstyle=\tiny\color{gray},
    stepnumber=1,
    showstringspaces=false,
    frame=single,
    breaklines=true,
    tabsize=4
}

\begin{document}

\renewcommand\Authfont{\large}
\renewcommand\Affilfont{\small}
\setlength{\affilsep}{0.5em}

\title{%
Central Limit Theorems for Sample Fréchet Means of Manifold-Valued Markov Chains
}

\author{Meshal Abuqrais}

\affil{%
Department of Mathematics, King's College London, London, United Kingdom\\
Department of Mathematics, Kuwait University, Kuwait\\
\vspace{0.5em}
\textit{Emails:}
\href{mailto:meshal.abuqrais@kcl.ac.uk}{\texttt{meshal.abuqrais@kcl.ac.uk}},
\href{mailto:m.abuqrais@ku.edu.kw}{\texttt{m.abuqrais@ku.edu.kw}}
}

\date{}

\maketitle

\begin{abstract}
In this article, we establish central limit theorems for sample Fréchet means of stationary ergodic Markov chains taking values in manifolds, extending the asymptotic theory previously developed for independent observations to a class of dependent manifold-valued processes. 
Our results derive the asymptotic normality of the sample Fréchet mean from a central limit condition at the population Fréchet mean, under suitable local regularity conditions.
We further provide sufficient geometric and probabilistic conditions under which these assumptions hold, formulated in terms of curvature bounds and a Wasserstein mixing condition.
As an application, we establish a central limit theorem for sample Fréchet means for a class of random dynamical systems generated by contractive random maps.
\end{abstract}

\section{Introduction}
The analysis of random variables on nonlinear spaces has become increasingly important in applications where observations possess an intrinsic geometric structure. 
Directional data, rotations, shapes, positive-definite matrices, and other structured objects are naturally represented as points in metric spaces or Riemannian manifolds rather than as vectors in Euclidean space; see, for example, \citep{mardia2009directional}, \citep{patrangenaru2016nonparametric}, and \citep{pennec2019riemannian}. 
In these settings, standard Euclidean statistical procedures are generally not intrinsic and may fail to respect the geometry of the underlying space. 
The Fréchet mean provides a natural generalisation of the Euclidean mean by defining centrality through the minimisation of the expected squared distance.

In the seminal works \citep{bhattacharya2003large} and \citep{bhattacharya2005large}, Bhattacharya and Patrangenaru studied fundamental asymptotic properties of sample Fréchet means.
They established strong consistency for random variables taking values in metric spaces under suitable conditions, as well as a central limit theorem for manifold-valued random variables.
These results were obtained under the assumption that the observations are independent and identically distributed.
A further central limit theorem in the independent and identically distributed setting was established in \citep{hotz2024central} for compact manifolds using a different proof strategy.
A central limit theorem for independent but not necessarily identically distributed manifold-valued random variables was established in \citep{kendall2011limit}.

To the best of our knowledge, no central limit theorem for sample Fréchet means under Markov dependence has yet been established.
Progress in this direction was made in \citep{abuqrais2026autoregressive}, where a strong law of large numbers was proved for sample Fréchet means of Markov-dependent random variables taking values in proper metric spaces.
More general strong laws for sample Fréchet means were established in \citep{jaffe2024fr}. 
These results naturally motivate the question of whether a corresponding central limit theorem can be obtained under Markov dependence.

The purpose of this paper is to establish central limit theorems for empirical Fréchet means associated with stationary ergodic manifold-valued Markov chains.
We prove two central limit theorems. The first assumes a central limit condition at the population Fréchet mean, formulated in a smooth local coordinate chart, while the second provides sufficient conditions for such a central limit theorem based on Wasserstein mixing rates and geometric assumptions.
The proof of the first theorem follows the strategy used by Bhattacharya and Patrangenaru in the independent and identically distributed setting \citep{bhattacharya2005large}.
For the second theorem, based on Wasserstein mixing rates, we establish intermediate lemmas that separate the contributions of the geometry of the underlying manifold and the probabilistic dependence structure. 
The geometric conditions concern the regularity of the exponential and logarithm maps and the curvature properties of the manifold, whereas the probabilistic conditions are expressed through the mixing rate of the Markov chain.
This allows the two components to be analysed separately before being combined in the proof of the central limit theorem. 
As an application of the general theory, we consider a class of Markov chains arising from random dynamical systems.

The paper is organised as follows. 
We first introduce the required geometric and probabilistic background. We then develop the asymptotic theory of empirical Fréchet means for stationary ergodic Markov chains and provide sufficient conditions based on Wasserstein mixing rates and curvature bounds. 
Finally, we apply the general results to random dynamical systems generated by contractive random maps.
\section{Preliminaries}

\subsection{Riemannian manifolds}
In this subsection, we introduce the core concepts of Riemannian geometry used in our work.
All manifolds considered herein are assumed to be finite-dimensional. 
For a comprehensive treatment of the subject, the reader is directed to \citep{Rie-Lee,jost2008riemannian}.

Let $M$ denote an $m$-dimensional smooth manifold. 
At any given point $p \in M$, we can attach a $m$-dimensional real vector space known as the tangent space, denoted by $T_pM$. 
Intuitively, this space comprises the velocity vectors of all possible smooth curves passing through $p$. 
We define the tangent bundle, $TM$, as the disjoint union of all such tangent spaces across the manifold. 
The bundle $TM$ possesses a canonical smooth structure such that the base projection map $\pi: TM \to M$, given by $(v; p) \mapsto \pi(v; p) = p$, is everywhere smooth.

A Riemannian metric on $M$ is defined as a smooth assignment of a positive-definite inner product, denoted $\inpr{\cdot}{\cdot}_p$, to the tangent space $T_pM$ at each point $p$. 
A manifold equipped with such a metric is called a Riemannian manifold. 
This inner product naturally induces a norm on each tangent space, written as $\norm{\cdot}_p$ or simply $\norm{\cdot}$ when there is no ambiguity.

If $M$ is connected, the Riemannian metric endows the manifold with a standard metric space structure. 
Given a piecewise smooth curve $\gamma:[0,1]\to M$, its length is defined as
\begin{equation*}
L(\gamma) = \int_0^1 \norm{\gamma'(t)}_{\gamma(t)} dt.
\end{equation*}
For $p,q\in M$, the distance function is defined by $d(p,q) = \inf_{\gamma} L(\gamma)$, where the infimum is taken over all piecewise smooth paths from $p$ to $q$. 
When referring to the Riemannian manifold as a metric space, it is with respect to this induced distance function.

We now turn to the construction of the exponential map, a fundamental tool that provides $M$ with natural local coordinates.
To proceed, we define geodesics.
Geodesics can be defined as the critical points of the energy functional 
\begin{equation*}
E(\gamma) = \frac{1}{2} \int_0^1 \norm{\gamma'(t)}_{\gamma(t)}^2 dt,
\end{equation*}
over the space of piecewise smooth curves with fixed endpoints. 
For any fixed point $p \in M$ and tangent vector $v \in T_pM$, there exists a unique maximal geodesic, denoted $\gamma_v(t; p)$, defined on an interval in $\mathbb{R}$ containing zero, subject to the initial conditions $\gamma_v(0; p) = p$ and $\gamma'_v(0; p) = v$. 
Assuming this geodesic is defined up to time $1$, the exponential map at $p$ is evaluated as
\begin{equation*}
\exp_p(v) = \gamma_v(1; p).
\end{equation*}

In general, the exponential map may not be globally defined over the entire tangent bundle $TM$. 
Nevertheless, for every $p$, $\exp_p$ is a local diffeomorphism on a neighbourhood of $0 \in T_pM$. 
The maximal radius of an open ball centred at $0$ in $T_pM$ on which the exponential is a diffeomorphism is called the injectivity radius of $M$ at $p$, written as $\injr{p}$.
The global injectivity radius of the manifold is correspondingly defined as $\injr{M} = \inf_{p \in M} \injr{p}$. 
Within the injectivity radius, the exponential map is invertible, and its inverse is called the logarithm map, denoted by $\log_p(\cdot)$.
Moreover, within the injectivity radius, the exponential map is a radial isometry.
In other words, for all $q \in B_p(r)$ with $r < \injr{p}$,
\begin{equation*}
d(p,q) = \norm{\log_p q}.
\end{equation*}
A ball $B_{p}(r)$ with $r\leq \injr{p}$ is called a geodesic ball centred at $p$.

A fundamental result in Riemannian geometry is the Hopf–Rinow theorem \citep{jost2008riemannian}, which guarantees that the metric completeness of $M$ is equivalent to geodesic completeness, that is, the exponential map is globally defined over the entire tangent bundle $TM$. 
Furthermore, these conditions are equivalent to $M$ satisfying the Heine-Borel property, meaning that a subset of $M$ is compact if and only if it is both closed and bounded.

Finally, we introduce the Riemannian curvature tensor and the sectional curvature.
By the fundamental theorem of Riemannian geometry, the metric uniquely determines a torsion-free and metric-compatible affine connection $\nabla$, known as the Levi-Civita connection. 
Let $\mathfrak{X}(M)$ denote the space of smooth vector fields on $M$. 
Using this connection, the Riemannian curvature tensor $R(X,Y): \mathfrak{X}(M) \to \mathfrak{X}(M)$ for any $X, Y \in \mathfrak{X}(M)$ is defined by
\begin{equation*}
R(X,Y)Z = \nabla_X \nabla_Y Z - \nabla_Y \nabla_X Z - \nabla_{[X,Y]} Z.
\end{equation*}
For a given point $p \in M$ and a two-dimensional plane spanned by linearly independent tangent vectors $X, Y \in T_pM$, the sectional curvature at $p$ is given by
\begin{equation*}
K_p(X,Y) = \frac{\inpr{R(X,Y)Y}{X}}{\norm{X}^2 \norm{Y}^2 - \inpr{X}{Y}^2}.
\end{equation*}
Generally, the sectional curvature varies depending on both the base point $p$ and the specific plane spanned by $X$ and $Y$. 
We say that $M$ has its curvature bounded above by a constant $K_U$ if $K_p(X,Y) \leq K_U$ for all points $p \in M$ and all pairs of linearly independent vectors $X, Y \in T_pM$. 
An analogous definition holds for curvature bounded from below. 
If the sectional curvature is constant across all points $p$ and all tangent planes, the manifold is said to have constant curvature $K$. 
The prototypical examples of spaces with constant curvature, known as space forms, are the flat Euclidean space $\mathbb{R}^m$ ($K = 0$), the unit sphere $\mathbb{S}^m$ equipped with the standard round metric ($K = 1$), and the hyperbolic space $\mathbb{H}^m$ endowed with the canonical hyperbolic metric ($K = -1$).

\subsection{Markov chains and Fréchet functions}
In this subsection, we introduce Markov chains, Wasserstein spaces, and the Fréchet function in metric spaces.
For comprehensive treatments of Markov chains and the Fréchet function of random objects, we refer the reader to \citep{benaim2022markov} and \citep{patrangenaru2016nonparametric}, respectively.
For Wasserstein spaces, we refer to \citep{villani2009optimal} or \citep{figalli2021invitation}.
Throughout this section, we assume that $(M,d)$ is a Polish metric space, i.e., complete and separable, equipped with its Borel $\sigma$-algebra $\mathcal{B}(M)$.
We also fix a filtered probability space $(\Omega,\mathcal{A},(\mathcal{A}_n)_{n\in\mathbb{N}_0},P)$.

We begin by defining Markov kernels on $M$.
Let $\mathsf{Prob}(M)$ denote the set of all probability measures on $M$.
A Markov kernel $\mathcal{K}$ on $(M,\mathcal{B}(M))$ is a map
\begin{equation*}
\mathcal{K}\colon M\times \mathcal{B}(M)\to [0,1],
\end{equation*}
such that, for all $B\in\mathcal{B}(M)$, the map $x \mapsto \mathcal{K}(x,B)$ is measurable and, for all $x\in M$, $\mathcal{K}(x,\cdot)\in \mathsf{Prob}(M)$.
An adapted stochastic process $(X_n)_{n\geq 0}$ of $M$-valued random variables on $(\Omega,\mathcal{A},(\mathcal{A}_n)_{n\geq 0},P)$ is called a Markov chain on $M$ with kernel $\mathcal{K}$ if for all $n\in\mathbb{N}_0$ and $B\in\mathcal{B}(M)$,
\begin{equation*}
P(X_{n+1}\in B \mid \mathcal{A}_n) = \mathcal{K}(X_n,B).
\end{equation*}
A Markov kernel $\mathcal{K}$ acts on $\mathsf{Prob}(M)$ by
\begin{equation*}
\mu \mapsto \mu\mathcal{K}(\cdot)=\int_{M}\mathcal{K}(x,\cdot)\,d\mu(x),
\end{equation*}
for all $\mu\in\mathsf{Prob}(M)$. 
Fixed points of this map are called invariant measures for $\mathcal{K}$ or the Markov chain.

Suppose that $\lambda$ is an invariant measure of the Markov chain $(X_n)_{n\geq 0}$.
The Markov chain is said to be stationary (or $\lambda$-stationary) if $\mathcal{L}(X_0) = \lambda$.
An invariant measure $\lambda$ of a Markov kernel $\mathcal{K}$ is said to be ergodic if for every bounded and measurable function $f$ on $M$ satisfying
\begin{equation*}
\int_{M}f(y)\mathcal{K}(x,dy)=f(x),
\end{equation*}
$\lambda$-a.s., it holds that $f$ is $\lambda$-a.s. constant.

To rigorously study the convergence of probability measures with respect to the underlying space, we equip subsets of $\mathsf{Prob}(M)$ with a metric space structure induced by the Wasserstein metric.
In order to define it, we need the concept of a coupling of measures.
Let $\mu,\nu \in \mathsf{Prob}(M)$ be two probability measures.
A coupling $\gamma$ of $\mu$ and $\nu$ is a probability measure on $M\times M$, i.e., $\gamma\in \mathsf{Prob}\left(M\times M\right)$, with the property that for all $A, B \in \mathcal{B}(M)$
\begin{equation*}
\gamma(A\times M)=\mu(A)\text{ and }\gamma(M\times B)=\nu(B).
\end{equation*}
The set of all couplings between $\mu$ and $\nu$ is denoted by $\Pi(\mu,\nu)$.

For $p\geq 1$, we define the space of probability measures with finite $p$th moments as
\begin{equation*}
\mathsf{Prob}_p(M)=\set{\mu\in \mathsf{Prob}(M)}{\int_{M}d(x,x_{0})^p d\mu(x)<\infty, \text{ for some $x_0\in M$}}.
\end{equation*}
The $p$-Wasserstein distance is a distance function on $\mathsf{Prob}_p(M)$ defined by
\begin{equation*}\label{Wasserstein distance}
\mathcal{W}_{p}(\mu,\nu)=\left(\inf_{\gamma\in\Pi(\mu,\nu)}\int_{M\times M}d(x,y)^{p}d\gamma(x,y)\right)^{1/p},
\end{equation*}
for all $\mu$ and $\nu$ in $\mathsf{Prob}_p(M)$.
We are primarily interested in the cases of $p=1,2$.
The Kantorovich-Rubinstein formula states that if $X$ and $Y$ are random variables whose probability measures are $\mu$ and $\nu$, then
\begin{equation*}
\mathcal{W}_1(\mu,\nu)=\inf_{\gamma\in \Pi(\mu,\nu)}E_\gamma[d(X,Y)] =\sup_{f}\left\{\int_{M}fd\mu - \int_{M}f d\nu\right\},
\end{equation*}
where the supremum is taken over all $\mathbb{R}$-valued $1$-Lipschitz functions $f$ on $M$.
It is worth mentioning that, since $M$ is a Polish space, the Wasserstein space $(\mathsf{Prob}_p(M),\mathcal{W}_p)$ is also Polish for all $1\leq p<\infty$; see Theorem 6.18 in \citep{villani2009optimal}.

Having established the metric framework for probability measures, we now turn to introducing statistical measures for $M$-valued random variables.
In a general metric space, there is no natural algebraic operation corresponding to the addition of points. 
Nevertheless, classical statistical quantities such as the mean must be defined through their variational characterisations.
If $X$ is an $M$-valued random variable, its Fréchet $r$th moment function is defined by
\begin{equation*}
\mathcal{F}_{r,X}(p) = \int_{M} d^r(x, p) P_X(dx) = E[d^r(X, p)],
\end{equation*}
where $r > 0$ and $P_X$ is the law of $X$.
The case of $r = 2$ is of special interest, and so we write $\mathcal{F}_X(p)$ or simply $\mathcal{F}(p)$ when there is no confusion.
We define the Fréchet mean set of an $M$-valued random object $X$ as
\begin{equation*}
\mathcal{U}_X= \operatorname{arg\,min}_{p \in M} \mathcal{F}(p).
\end{equation*}
It can be shown that this definition of the mean set coincides with the notion of mean in linear spaces.
However, it may be the case that there is more than one mean point, which does not occur in Euclidean spaces.
Whether the mean is unique or not depends on the structure of the space.
For instance, for Riemannian manifolds uniqueness may be guaranteed if the random variable is supported in a sufficiently small geodesic ball; see \citep{afsari2011riemannian}.
If the space is a Hadamard manifold, global uniqueness holds, see, for
example, \citep{pennec2006intrinsic} or \citep{sturm2003probability}.

The natural estimator for the Fréchet mean set based on $X_1, \ldots, X_n$ is
\begin{equation*}
\mathcal{U}_n = \operatorname{arg\,min}_{p \in M} \frac{1}{n} \sum_{k=1}^n d^2(X_k, p),
\end{equation*}
which is the Fréchet mean set of the empirical Fréchet function $\mathcal{F}_n(\cdot)=\frac{1}{n}\sum_{k=1}^n d^2(X_k, \cdot)$.
It was shown in \citep{bhattacharya2003large} that, under the assumption that $X_1, \ldots, X_n$ are independent and identically distributed, together with suitable assumptions on the space, the Fréchet sample mean set is strongly consistent (i.e., converges almost surely) for the true Fréchet mean set under an appropriate notion of set convergence.
Such a strong law of large numbers for the sample Fréchet mean holds in a more general setting.
Specifically, the same law holds under Markov dependence; see \citep{abuqrais2026autoregressive} and \citep{jaffe2024fr}.
This strong law for the sample Fréchet mean for Markov chains will play a key role in the proof of our central limit theorems.
Throughout this section, we assume that $(\Omega,\mathcal{A},P)$ is the underlying probability space on which the Markov chains are defined.

\section{Central Limit Theorems}
The primary objective of this section is to establish two central limit theorems for the sample Fréchet mean of a Markov chain on a manifold. 
We assume throughout that the population Fréchet mean is unique.
The first theorem assumes a central limit condition at the Fréchet mean, formulated in a local coordinate chart.
Since this condition may be difficult to verify directly, the second theorem provides sufficient conditions in terms of mixing rates and curvature bounds.
Let $V$ be a finite $m$-dimensional real vector space and let $f$ be a smooth $\mathbb{R}$-valued map on $V$.
By $Df$ we mean the differential of $f$.
That is, $Df(v)(\cdot)$ gives the directional derivative of $f$ at $v$.
More generally, if $F:V\to N$ is a smooth map into a smooth manifold $N$, we write $DF|_v$ for its differential at $v$.
Similarly, by $D^2f$ we mean the Euclidean Hessian operator of $f$.
That is, $D^2f$ sends $v$ to the Hessian operator
\begin{equation*}
D^2f(v):V\times V\longrightarrow\mathbb{R},
\end{equation*}
which sends $(u,w)$ to $D^2f(v)(u,w)$ for all $u,w\in V$.
As a bilinear form, its operator norm is defined by
\begin{equation*}
\norm{D^2f(v)}_{\mathrm{op}}=\sup_{\substack{\norm{u}\leq 1,\ \norm{w}\leq 1}}
\abs{D^2f(v)(u,w)}.
\end{equation*}
If $f$ is a function of two variables, $v_1$ and $v_2$, we write $Df(v_1;v_2)$ to indicate that we are differentiating with respect to the first variable $v_1$ while holding $v_2$ fixed.
For a smooth $\mathbb{R}$-valued function $h$ on a Riemannian manifold $M$, we write $\operatorname{Hess}h$ for its Riemannian Hessian.
For a smooth manifold-valued map $F:V\to M$, we write $\nabla DF$ for the covariant derivative of its differential with respect to the Levi-Civita connection on $M$.
\subsection{Central limit theorem in a smooth chart}
Our first central limit theorem is formulated in terms of a local central limit condition in a smooth coordinate chart.
In particular, the required regularity assumptions are imposed only in a neighbourhood of the population Fréchet mean.
The proof follows closely the argument of \citep{bhattacharya2005large}.

\begin{theorem}\label{thm:local clt for sample Frechet mean}
Let $(M,d)$ be a proper metric space that is also a $m$-dimensional smooth manifold, and assume that the metric topology agrees with the manifold topology.
Let $(X_n)_{n\geq 0}$ be an $M$-valued Markov chain and let $\lambda$ be an ergodic invariant measure of the chain.
Assume that the chain is stationary with respect to $\lambda$ ($X_0\sim\lambda$).
Suppose that $\lambda$ has finite Fréchet function $\mathcal F$ and has $\mu\in M$ as its unique Fréchet mean.
Let $\mathcal F_n$ be the empirical Fréchet function and $\mu_n\in \operatorname{arg\min}_{p\in M}\mathcal F_n(p)$.
For a neighbourhood $U$ of $\mu$, let $\phi:U\to\phi(U)\subset\mathbb{R}^m$ be a smooth coordinate chart centred at $\mu$, $\phi(\mu)=0$.
Fix $r>0$ such that
\begin{equation*}
\{v\in\mathbb{R}^m:\norm{v}<r\}
\subset \phi(U).
\end{equation*}
Define, for $\norm{v}<r$,
\begin{equation*}
F(v)=\mathcal F(\phi^{-1}(v)),
\qquad
F_n(v)=\mathcal F_n(\phi^{-1}(v)).
\end{equation*}
Let
\begin{equation*}
G(v;x)=d^2(\phi^{-1}(v),x).
\end{equation*}
Assume that, for $\lambda$-almost every $x\in M$, the map
\begin{equation*}
v\longmapsto G(v;x)
\end{equation*}
is $C^2$ on $\{v\in\mathbb{R}^m:\norm{v}<r\}$.
Assume that the Hessian $D^2F(0)$ is invertible and that, for all $0<\delta<r$,
\begin{equation}\label{Integrability tail-curvature condition}
E_\lambda\left[\sup_{\norm{v}<\delta}\norm{D^2G(v;X_0)}_{\mathrm{op}}\right]<\infty.
\end{equation}
Assume that
\begin{equation}\label{CLT at one point}
\frac{1}{\sqrt n}\sum_{k=1}^n DG(0;X_k)\longrightarrow N_m(0,\Sigma),
\end{equation}
in distribution in $\mathbb{R}^m$, where $\Sigma$ is a covariance matrix.
Then, as $n\to\infty$,
\begin{equation*}
\sqrt n\,\phi(\mu_n) \longrightarrow N_m\left(0,(D^2F(0))^{-1}\Sigma(D^2F(0))^{-1}\right),
\end{equation*}
in distribution.
\end{theorem}
Before we prove the theorem, we need the following lemma which ensures the convergence of the Hessians.
\begin{lemma}
\label{lemma: Hessian a.s. uniform convergence}
Let $(M,d)$ be a proper metric space that is also a $m$-dimensional smooth manifold, and assume that the metric topology agrees with the manifold topology.
Fix $\mu\in M$.
Let $\{X_n\}_{n\geq 0}$ be an $M$-valued Markov chain with an ergodic invariant measure $\lambda$ whose Fréchet function is finite everywhere, and assume that $\{X_n\}_{n\geq 0}$ is stationary with respect to $\lambda$.
Let $\phi:U\to\phi(U)\subset\mathbb{R}^m$ be a smooth coordinate chart containing $\mu$, with $\phi(\mu)=0$.
Define
\begin{equation*}
F_n(v)=\frac{1}{n}\sum_{k=1}^nd^2(\phi^{-1}(v),X_k),
\qquad
F(v)=E_\lambda[d^2(\phi^{-1}(v),X_0)].
\end{equation*}
Let
\begin{equation*}
G(v;x)=d^2(\phi^{-1}(v),x).
\end{equation*}
Fix $\delta>0$, and let
\begin{equation*}
K_\delta=\{v\in\mathbb{R}^m:\norm{v}\leq\delta\}
\subset\phi(U).
\end{equation*}
Assume that, for $\lambda$-almost every $x\in M$, the map $v\longmapsto G(v;x)$ is $C^2$ on an open neighbourhood of $K_\delta$, and that
\begin{equation*}
E_\lambda\left[\sup_{v\in K_\delta}\norm{D^2G(v;X_0)}_{\mathrm{op}}\right]<\infty.
\end{equation*}
Then
\begin{equation*}
\sup_{\norm{v}\leq\delta}
\norm{D^2F_n(v)-D^2F(v)}_{\mathrm{op}}
\overset{P\text{-a.s.}}{\longrightarrow}
0,
\end{equation*}
as $n\to\infty$.
\end{lemma}
\begin{proof}
Fix $\varepsilon>0$.
For any $v,w\in K_\delta$, we have
\begin{equation*}
\begin{aligned}
\norm{D^2F_n(v)-D^2F(v)}_{\mathrm{op}}
&\leq
\norm{D^2F_n(v)-D^2F_n(w)}_{\mathrm{op}}
\\
&\quad+
\norm{D^2F(w)-D^2F(v)}_{\mathrm{op}}
\\
&\quad+
\norm{D^2F_n(w)-D^2F(w)}_{\mathrm{op}}.
\end{aligned}
\end{equation*}
We estimate the three terms on the right-hand side separately to obtain a uniform bound for all $v\in K_\delta$.
First, define the function
\begin{equation*}
H(x)=\sup_{v\in K_\delta}\norm{D^2G(v;x)}_{\mathrm{op}},
\end{equation*}
which satisfies $H\in L^1(\lambda)$ by assumption.
For $t>0$, define the function
\begin{equation*}
A_\delta(t;x)=\sup_{\substack{v,w\in K_\delta\\ \norm{v-w}\leq t}}
\norm{D^2G(v;x)-D^2G(w;x)}_{\mathrm{op}}.
\end{equation*}
Since, for $\lambda$-almost all $x\in M$, the map $v\mapsto G(v;x)$ is $C^2$ on $K_\delta$, the map $v\mapsto D^2G(v;x)$ is continuous.
Thus,
\begin{equation*}
\lim_{t\to 0}A_\delta(t;x)=0
\end{equation*}
for $\lambda$-almost all $x\in M$.
Since $A_\delta(t;x)\leq 2H(x)$ and $H\in L^1(\lambda)$, the Dominated Convergence Theorem implies
\begin{equation*}
\lim_{t\to 0}E_\lambda[A_\delta(t;X_0)]=0.
\end{equation*}
Therefore, there exists $\eta>0$ such that
\begin{equation*}
E_\lambda[A_\delta(\eta;X_0)]<\frac{\varepsilon}{3}.
\end{equation*}
For the first term, for any $v,w\in K_\delta$ such that $\norm{v-w}\leq\eta$,
\begin{equation*}
\norm{D^2F_n(v)-D^2F_n(w)}_{\mathrm{op}}\leq\frac{1}{n}\sum_{k=1}^n A_\delta(\eta;X_k).
\end{equation*}
Because $A_\delta(\eta;\cdot)\in L^1(\lambda)$, the Ergodic Theorem implies that there exists a set $\Omega_A$ with $P(\Omega_A)=1$ such that for every $\omega\in\Omega_A$,
\begin{equation*}
\lim_{n\to\infty}\frac{1}{n}\sum_{k=1}^n A_\delta(\eta;X_k(\omega))=E_\lambda[A_\delta(\eta;X_0)]
<\frac{\varepsilon}{3}.
\end{equation*}
Next, we bound the second term.
Since
\begin{equation*}
\norm{D^2G(v;X_0)}_{\mathrm{op}}
\leq
H(X_0)
\end{equation*}
and $H\in L^1(\lambda)$, the Dominated Convergence Theorem yields
\begin{equation*}
D^2F(v)=E_\lambda[D^2G(v;X_0)].
\end{equation*}
Furthermore, by Jensen's inequality, for any $v,w\in K_\delta$ with $\norm{v-w}\leq\eta$,
\begin{equation*}
\norm{D^2F(v)-D^2F(w)}_{\mathrm{op}}\leq E_\lambda[A_\delta(\eta;X_0)]<\frac{\varepsilon}{3}.
\end{equation*}
Now we bound the last term.
By the compactness of $K_\delta$, there exist finitely many points $w_1,\ldots,w_J\in K_\delta$ such that
\begin{equation*}
K_\delta \subset\bigcup_{j=1}^J B(w_j,\eta).
\end{equation*}
The bound
\begin{equation*}
\norm{D^2G(w_j;X_0)}_{\mathrm{op}}\leq H(X_0)
\end{equation*}
guarantees that $D^2G(w_j;\cdot)\in L^1(\lambda)$ for each centre $w_j$.
By the Ergodic Theorem, for each $j=1,\ldots,J$, there exists $\Omega_j$ with $P(\Omega_j)=1$ such that for all $\omega\in\Omega_j$,
\begin{equation*}
\lim_{n\to\infty}D^2F_n(w_j)(\omega)=D^2F(w_j).
\end{equation*}
Define
\begin{equation*}
\Omega_*=\Omega_A\cap\bigcap_{j=1}^{J}\Omega_j.
\end{equation*}
Since this is a finite intersection of full-measure sets, $P(\Omega_*)=1$.
Let $\omega\in\Omega_*$.
There exists $N_1(\omega)$ such that for all $n\geq N_1(\omega)$,
\begin{equation*}
\max_{1\leq j\leq J}\norm{D^2F_n(w_j)(\omega)-D^2F(w_j)}_{\mathrm{op}}<\frac{\varepsilon}{3}.
\end{equation*}
There also exists $N_2(\omega)$ such that for all $n\geq N_2(\omega)$,
\begin{equation*}
\frac{1}{n}\sum_{k=1}^n A_\delta(\eta;X_k(\omega))<\frac{\varepsilon}{3}.
\end{equation*}
Let
\begin{equation*}
N(\omega)=\max\{N_1(\omega),N_2(\omega)\}.
\end{equation*}
For any $v\in K_\delta$, there exists $j\in\{1,\ldots,J\}$ such that $\norm{v-w_j}<\eta$.
For all $n\geq N(\omega)$, applying the uniform bounds yields
\begin{equation*}
\begin{aligned}
\norm{D^2F_n(v)-D^2F(v)}_{\mathrm{op}}
&\leq
\norm{D^2F_n(v)-D^2F_n(w_j)}_{\mathrm{op}}
\\
&\quad+
\norm{D^2F(w_j)-D^2F(v)}_{\mathrm{op}}
\\
&\quad+
\norm{D^2F_n(w_j)-D^2F(w_j)}_{\mathrm{op}}
\\
&<\frac{\varepsilon}{3}+\frac{\varepsilon}{3}+\frac{\varepsilon}{3}=\varepsilon.
\end{aligned}
\end{equation*}
Therefore,
\begin{equation*}
\sup_{v\in K_\delta}\norm{D^2F_n(v)-D^2F(v)}_{\mathrm{op}}\to 0
\end{equation*}
$P$-a.s. as $n\to\infty$.
\end{proof}

\begin{proof}[Proof of Theorem \ref{thm:local clt for sample Frechet mean}]
By the assumed CLT,
\begin{equation*}
\frac{1}{\sqrt n}\sum_{k=1}^nDG(0;X_k)\xrightarrow{d}N_m(0,\Sigma).
\end{equation*}
By the SLLN for the sample Fréchet mean \citep{abuqrais2026autoregressive},
$\mu_n\to\mu$ $P$-a.s.
Therefore, for $P$-almost every $\omega$, $\mu_n(\omega)\in U$ for all sufficiently large $n$.
For all such $n$, set $v_n=\phi(\mu_n).$
Then, $v_n\to0$ $P$-a.s.
Note that $DF_n(v_n)=0$ since $\mu_n$ is an empirical Fréchet minimiser.
For each such sufficiently large $n$, consider the path $\tilde{\gamma}_n:[0,1]\to\mathbb{R}^m$ given by $\tilde{\gamma}_n(t)=tv_n$, and define
\begin{equation*}
\Phi_n(t)=DF_n(\tilde{\gamma}_n(t)).
\end{equation*}
Then $\Phi_n(1)=0$ and
\begin{equation*}
\Phi_n(0)=DF_n(0)=\frac{1}{n}\sum_{k=1}^n DG(0;X_k).
\end{equation*}
Moreover,
\begin{align*}
\Phi_n(1)-\Phi_n(0)&=\int_0^1\Phi_n'(t)\,dt=\int_0^1 D^2F_n\big|_{\tilde{\gamma}_n(t)}[\tilde{\gamma}_n'(t)]\,dt
\\
&=\left(\int_0^1D^2F_n\big|_{tv_n}dt\right)[v_n].
\end{align*}
Hence,
\begin{equation*}
-\frac{1}{n}\sum_{k=1}^nDG(0;X_k)=\left(\int_0^1D^2F_n\big|_{tv_n}\,dt\right)[v_n].
\end{equation*}
Set
\begin{equation*}
A_n=\int_0^1D^2F_n\big|_{tv_n}dt.
\end{equation*}
Then
\begin{equation*}
-\frac{1}{n}\sum_{k=1}^n DG(0;X_k)=A_n[v_n].
\end{equation*}
Fix $0<\delta<r$.
Since $v_n\to0$ $P$-a.s., for $P$-almost every $\omega$, there exists $N$ such that for all $n\geq N$ and all $t\in[0,1]$, $\norm{tv_n}<\delta$.
For each such $n$ and $t$, the triangle inequality gives
\begin{equation*}
\norm{D^2F_n|_{tv_n}-D^2F|_0}_{\mathrm{op}}
\leq
\norm{D^2F_n|_{tv_n}-D^2F|_{tv_n}}_{\mathrm{op}}
+
\norm{D^2F|_{tv_n}-D^2F|_0}_{\mathrm{op}}.
\end{equation*}
Taking the supremum over $t\in[0,1]$ gives
\begin{align*}
\sup_{t\in[0,1]}
\norm{D^2F_n|_{tv_n}-D^2F|_0}_{\mathrm{op}}
&\leq
\sup_{\norm{v}\leq\delta}
\norm{D^2F_n|_v-D^2F|_v}_{\mathrm{op}}
\\
&\quad+
\sup_{t\in[0,1]}
\norm{D^2F|_{tv_n}-D^2F|_0}_{\mathrm{op}}.
\end{align*}
By Lemma \ref{lemma: Hessian a.s. uniform convergence} and $v_n\to0$ $P$-a.s.,
\begin{equation*}
\sup_{t\in[0,1]}\norm{D^2F_n\big|_{tv_n}-D^2F\big|_0}_{\mathrm{op}}\overset{P\text{-a.s.}}{\longrightarrow}0.
\end{equation*}
Therefore,
\begin{equation*}
A_n\overset{P\text{-a.s.}}{\longrightarrow}D^2F(0).
\end{equation*}
Since $D^2F(0)$ is invertible, $A_n$ is invertible for all sufficiently large $n$, $P$-a.s., and
\begin{equation*}
A_n^{-1}\overset{P\text{-a.s.}}{\longrightarrow}
\left(D^2F(0)\right)^{-1}.
\end{equation*}
Moreover,
\begin{equation*}
\sqrt n\,v_n=-A_n^{-1}\left(\frac{1}{\sqrt n}\sum_{k=1}^nDG(0;X_k)\right).
\end{equation*}
Therefore, by Slutsky's theorem,
\begin{equation*}
\sqrt n\,\phi(\mu_n)\longrightarrow N_m\left(0,\left(D^2F(0)\right)^{-1}\Sigma\left(D^2F(0)\right)^{-1}\right),
\end{equation*}
in distribution as $n\to\infty$.
\end{proof}

Theorem \ref{thm:local clt for sample Frechet mean} extends the framework of \citep{bhattacharya2005large} primarily by allowing the observations to form a stationary ergodic Markov chain, rather than requiring them to be i.i.d. 
In addition, we localise their regularity condition by requiring the squared-distance function to be $C^2$ only in a neighbourhood of the population Fréchet mean.
We now state the corresponding result for Riemannian manifolds as a corollary.

\begin{corollary}\label{corollary:clt for sample Frechet mean Riemannian}
Let $M$ be a complete $m$-dimensional Riemannian manifold, and suppose that the chain $\{X_n\}_{n\geq 0}$ is supported in $B_\mu(r)$, where $r\in(0,\infty]$, $B_\mu(r)$ is a geodesic ball [$B_\mu(\infty)=M$ if $r=\infty$].
Assume that the conditions of Theorem~\ref{thm:local clt for sample Frechet mean} hold with $U=B_\mu(r)$ and $\phi=\log_\mu$, with condition \eqref{CLT at one point} replaced by
\begin{equation*}
\frac{1}{\sqrt n}\sum_{k=1}^n \log_\mu X_k
\longrightarrow N_m(0,\Sigma)
\end{equation*}
in distribution.
Let
\begin{equation*}
F(v)=E_\lambda\!\left[d^2(\exp_\mu v,X_0)\right].
\end{equation*}
Then
\begin{equation*}
\sqrt n\,\log_\mu\mu_n
\longrightarrow
N_m\left(0,4(D^2F(0))^{-1}\Sigma(D^2F(0))^{-1}\right),
\end{equation*}
in distribution.
\end{corollary}

\begin{proof}
By the Hopf--Rinow theorem, a complete Riemannian manifold is a proper metric space.
Taking $\phi=\log_\mu$ on $B_\mu(r)$, so that $\phi^{-1}=\exp_\mu$, and identifying $T_\mu M$ with $\mathbb{R}^m$, we have
\begin{equation*}
DG(0;x)=-2\log_\mu x.
\end{equation*}
Hence, the assumed central limit theorem implies
\begin{equation*}
\frac{1}{\sqrt n}\sum_{k=1}^n DG(0;X_k)
\longrightarrow N_m(0,4\Sigma),
\end{equation*}
in distribution.
Thus condition \eqref{CLT at one point} of Theorem~\ref{thm:local clt for sample Frechet mean} holds with covariance matrix $4\Sigma$.
The conclusion therefore follows from Theorem~\ref{thm:local clt for sample Frechet mean}.
\end{proof}

\subsection{CLT based on a Wasserstein mixing rate}
The conditions \eqref{CLT at one point} and \eqref{Integrability tail-curvature condition} can be hard to verify in general.
We give sufficient conditions under which both hold.
Specifically, we show that the condition \eqref{Integrability tail-curvature condition} holds for manifolds with bounded curvature and that condition \eqref{CLT at one point} holds under a $\mathcal{W}_1$-mixing condition.
But first, we need a geometric lemma about the logarithm and exponential maps.
For $K\in\mathbb{R}$, define
\begin{equation*}
S_K(t)
=
\begin{cases}
\dfrac{\sin(\sqrt{K}\,t)}{\sqrt{K}}, & K>0, \\
t, & K=0,\\
\dfrac{\sinh(\sqrt{-K}\,t)}{\sqrt{-K}}, & K<0,
\end{cases}
\end{equation*}
and
\begin{equation*}
\operatorname{sinc}_K(t)
=
\begin{cases}
\dfrac{S_K(t)}{t}, & t>0,\\[0.5em]
1, & t=0
\end{cases}.
\end{equation*}

\begin{lemma}\label{lem:uniform-exp-log-bounds}
Let $M$ be a complete Riemannian manifold, and fix $\mu\in M$.
Assume that the curvature of $M$ is bounded above by $K_U$ and below by $K_L$ for some $K_L,K_U\in\mathbb{R}$.
Let $B_\mu(r)$ be the geodesic ball of radius
\begin{equation*}
r<\min\left\{\injr{\mu},\frac{\pi}{\sqrt{K_U}}\right\},
\end{equation*}
where $\pi/\sqrt{K_U}:=\infty$ when $K_U\leq0$.
Define
\begin{equation*}
\Delta_L=\max_{0\leq t\leq r}\operatorname{sinc}_{K_L}(t),
\qquad
\Delta_U= \min_{0\leq t\leq r}\operatorname{sinc}_{K_U}(t).
\end{equation*}
Then
\begin{equation*}
0<\Delta_U\leq1\leq\Delta_L<\infty.
\end{equation*}
Moreover, for every $v,w\in T_\mu M$ with $\lVert v\rVert<r$,
\begin{equation}\label{eq:dexp-two-sided-bound}
\Delta_U\lVert w\rVert
\leq
\left\lVert D\exp_\mu\big|_v(w)\right\rVert
\leq
\Delta_L\lVert w\rVert.
\end{equation}
Consequently,
\begin{equation*}
\left\lVert D\exp_\mu\big|_v\right\rVert_{\mathrm{op}}
\leq
\Delta_L,
\end{equation*}
and
\begin{equation*}
\left\lVert D\log_\mu\big|_p\right\rVert_{\mathrm{op}}
\leq
\frac{1}{\Delta_U},
\end{equation*}
where $p=\exp_{\mu}(v)$.
If $M$ is a Hadamard manifold, then we can take $r=\infty$ and $K_U=0$, in which case $\Delta_U=1$.
\end{lemma}

\begin{proof}
The assertion is immediate when $v=0$, since
\begin{equation*}
D\exp_\mu\big|_0=\operatorname{Id}_{T_\mu M}.
\end{equation*}
Suppose that $v\neq0$, and set $a=\lVert v\rVert$. 
Consider the unit-speed geodesic
\begin{equation*}
\gamma(u)=\exp_\mu\left(\frac{u}{a}v\right),
\qquad
0\leq u\leq a.
\end{equation*}
Decompose
\begin{equation*}
w=w^\top+w^\perp
\end{equation*}
into its components tangential and normal to $v$.
Let $J$ be the Jacobi field along $\gamma$ associated with the variation
\begin{equation*}
(s,u)\mapsto
\exp_\mu\left(\frac{u}{a}(v+s w)\right).
\end{equation*}
Then
\begin{equation*}
J(0)=0,
\qquad
D_uJ(0)=\frac{w}{a},
\qquad
J(a)=D\exp_\mu\big|_v(w).
\end{equation*}

The tangential component satisfies
\begin{equation*}
\left\lVert J^\top(a)\right\rVert
=
\left\lVert w^\top\right\rVert.
\end{equation*}
By the Rauch comparison theorem, the normal component satisfies
\begin{equation*}
\frac{S_{K_U}(a)}{a}
\left\lVert w^\perp\right\rVert
\leq
\left\lVert J^\perp(a)\right\rVert
\leq
\frac{S_{K_L}(a)}{a}
\left\lVert w^\perp\right\rVert.
\end{equation*}
Equivalently,
\begin{equation*}
\operatorname{sinc}_{K_U}(a)
\left\lVert w^\perp\right\rVert
\leq
\left\lVert J^\perp(a)\right\rVert
\leq
\operatorname{sinc}_{K_L}(a)
\left\lVert w^\perp\right\rVert.
\end{equation*}

Since the tangential and normal components remain orthogonal, it follows that
\begin{align*}
\norm{D\exp_\mu\big|_v(w)}^2
&=
\norm{J^\top(a)}^2
+
\norm{J^\perp(a)}^2\\
&\leq
\norm{w^\top}^2
+
\Delta_L^2\norm{w^\perp}^2\\
&\leq
\Delta_L^2\norm{w}^2.
\end{align*}
Similarly,
\begin{align*}
\norm{D\exp_\mu\big|_v(w)}^2
&\geq
\norm{w^\top}^2
+
\Delta_U^2\norm{w^\perp}^2\\
&\geq
\Delta_U^2\norm{w}^2.
\end{align*}
This proves \eqref{eq:dexp-two-sided-bound}.
For $p=\exp_\mu(v)$, the condition $\norm{v}<\injr{\mu}$ implies
\begin{equation*}
D\log_\mu\big|_p
=
\left(D\exp_\mu\big|_v\right)^{-1}.
\end{equation*}
The lower bound in \eqref{eq:dexp-two-sided-bound} therefore gives
\begin{equation*}
\norm{D\log_\mu\big|_p}_{\mathrm{op}}
\leq
\frac{1}{\Delta_U}.
\end{equation*}
In the Hadamard case, one may take $r=\infty$ and $K_U=0$ to obtain
\begin{equation*}
\Delta_U=\min_{0\leq t}\operatorname{sinc}_0(t)=1.
\end{equation*}
\end{proof}
In the Hadamard case, if the sectional curvature of $M$ is not bounded below,
then no global upper bound for the differential of the exponential map follows
from these assumptions. 
In this case, the argument yields only
\begin{equation*}
\norm{D\exp_\mu\big|_v(w)}\geq \norm{w},
\end{equation*}
for all $v,w\in T_\mu M$.

\begin{lemma}\label{lemma: bounded curvature and Hessian a.s. uniform convergence}
Let $M$ be a complete Riemannian manifold and fix $\mu\in M$.
Suppose that the curvature of $M$ is bounded above and below by $K_U$ and $K_L$, respectively, where $K_U,K_L\in\mathbb{R}$.
Let $B_\mu(r)$ be a geodesic ball of radius
\begin{equation*}
r<\frac{1}{2}\min\left\{\injr{M},\frac{\pi}{\sqrt{K_U}}\right\},
\end{equation*}
where $\frac{\pi}{\sqrt{K_U}}=\infty$, if $K_U\leq 0$. 
If $M$ is a Hadamard manifold, then we take $r=\infty$, in which case $B_{\mu}(\infty)=M$.
Let $(X_n)_{n\ge 0}$ be a Markov chain supported on $B_{\mu}(r)$, with an ergodic invariant measure $\lambda$ whose Fréchet function is finite everywhere, and assume that $(X_n)$ is stationary with respect to $\lambda$.
Define
\begin{equation*}
F_n(v)=\frac{1}{n}\sum_{k=1}^n d^2(\exp_\mu v,X_k),
\qquad
F(v)=E_\lambda[d^2(\exp_\mu v,X_0)].
\end{equation*}
Then, for every $0<\delta<r$,
\[
\sup_{\norm{v}\leq\delta}
\norm{D^2F_n(v)-D^2F(v)}_{\mathrm{op}}
\xrightarrow{\mathrm{a.s.}}0,
\]
as $n\to\infty$.
\end{lemma}

\begin{proof}
Define $G(v;x) = d^2(\exp_\mu v, x)$ and let $y = \exp_\mu v$. 
The first derivative of $G$ with respect to $v$ along a vector $W \in T_\mu M$ is
\begin{equation*}
DG(v; x)(W) = \Big\langle -2 \log_{y}x, D\exp_\mu|_v(W) \Big\rangle_{y}.
\end{equation*}
Differentiating again with respect to $v$ along a vector $Z \in T_\mu M$ yields the Hessian $D^2G(v;x)(Z, W)$. 
By the product rule for covariant derivatives, this splits into two terms
\begin{equation*}
\begin{aligned}
D^2G(v;x)(Z, W)
&= \operatorname{Hess}_y d^2(\cdot, x)
\left(D(\exp_\mu)|_v(Z), D\exp_\mu|_v(W) \right) \\
& + \inpr{-2 \log_{y}x}{
\nabla D\exp_\mu|_v(Z, W)}_{y}.
\end{aligned}
\end{equation*}

We bound the operator norm $\norm{D^2G(v;x)}_{\mathrm{op}}$ uniformly for all $v$ in the compact ball $K_\delta = \set{v \in T_\mu M}{\norm{v} \leq \delta}$. 
Let $\rho = d(y, x) = d(\exp_\mu v, x)$. 
By Lemma \ref{lem:uniform-exp-log-bounds}, since $\delta<r$, the first derivative of the exponential map is uniformly bounded by $\Delta_L(\delta)$, meaning
\begin{equation*}
\norm{D(\exp_\mu)|_v}_{\mathrm{op}} \leq \Delta_L(\delta),
\end{equation*}
where $\Delta_L(\delta)$ denotes the constant from Lemma \ref{lem:uniform-exp-log-bounds} corresponding to radius $\delta$.
Furthermore, since $\exp_\mu$ is smooth and $K_\delta$ is compact, the covariant derivative of its differential is uniformly bounded by some finite constant $B_\delta>0$
\begin{equation*}
\norm{\nabla D(\exp_\mu)|_v}_{\mathrm{op}} \le B_\delta.
\end{equation*}
Applying this bound $\Delta_L(\delta)$, the operator norm of the first term is bounded by $\norm{\operatorname{Hess}_y d^2(\cdot, x)}_{\mathrm{op}}\Delta_L^2(\delta)$.
Since
\begin{equation*}
\rho=d(\exp_\mu v,x)\leq d(\exp_\mu v,\mu)+d(\mu,x),
\end{equation*}
Theorem 6.6.1 of \citep{jost2008riemannian} implies that there exist deterministic constants $C_1,C_2>0$ such that
\begin{equation*}
\norm{\operatorname{Hess}_y d^2(\cdot,x)}_{\mathrm{op}}
\le C_1\rho+C_2.
\end{equation*}
Therefore, the first term is bounded by $(C_1\rho+C_2)\Delta_L^2(\delta)$.
The norm of the gradient vector $-2 \log_{y}x$ is exactly $2\rho$, so applying the bound $B_\delta$, the second term is bounded by $2\rho B_\delta$.

Combining these bounds, the operator norm of the full Hessian grows at most linearly in the distance $\rho$
\begin{equation*}
\norm{D^2G(v;x)}_{\mathrm{op}} \le C_3 \rho + C_4,
\end{equation*}
for some constants $C_3,C_4>0$.
By the triangle inequality, $\rho = d(\exp_\mu v, x) \le d(\exp_\mu v, \mu) + d(\mu, x)$.
Because $v \in K_\delta$, we have $d(\exp_\mu v, \mu) = \norm{v} \le \delta$.
Taking the supremum over $K_\delta$ gives
\begin{equation*}
\sup_{v \in K_\delta} \norm{D^2G(v; x)}_{\mathrm{op}} \le C_5 d(\mu, x) + C_6,
\end{equation*}
for some constants $C_5,C_6>0$.
Replacing the fixed point $x$ with $X_0 \sim \lambda$ and taking the expectation yields
\begin{equation*}
E_\lambda \left[ \sup_{v \in K_\delta} \norm{D^2G(v; X_0)}_{\mathrm{op}} \right] \le C_5 E_\lambda [ d(\mu, X_0) ] + C_6.
\end{equation*}
By hypothesis, the Fréchet function is finite everywhere,
\begin{equation*}
E_\lambda \left[ \sup_{v \in K_\delta} \norm{D^2G(v; X_0)}_{\mathrm{op}} \right] < \infty.
\end{equation*}
The conclusion now follows from Lemma~\ref{lemma: Hessian a.s. uniform convergence}.
\end{proof}

\begin{theorem}\label{thm:frechet-clt-w1-mixing}
Let $M$ be a complete $m$-dimensional Riemannian manifold and fix $\mu\in M$.
Suppose that the curvature of $M$ is bounded above and below by $K_U$ and $K_L$, respectively, where $K_U,K_L\in\mathbb{R}$.
Let $B_\mu(r)$ be the geodesic ball of radius
\begin{equation*}
r<\frac{1}{2}\min\left\{\injr{M},\frac{\pi}{\sqrt{K_U}}\right\},
\end{equation*}
where $\frac{\pi}{\sqrt{K_U}}=\infty$, if $K_U\leq 0$. 
If $M$ is a Hadamard manifold, then we take $r=\infty$, in which case $B_{\mu}(\infty)=M$.
Let $(X_n)_{n\ge 0}$ be a Markov chain supported on $B_{\mu}(r)$, with an ergodic invariant measure $\lambda$ whose Fréchet function is finite everywhere, and assume that $(X_n)$ is stationary with respect to $\lambda$.
Assume that $\mu$ is the unique Fréchet mean of $\lambda$ and let $\mu_n$ be the sample Fréchet mean.
Define, in a neighbourhood of $0\in T_\mu M$,
\begin{equation*}
F(v)=\mathcal F(\exp_\mu v),
\end{equation*}
and assume $D^2F(0)$ is invertible.
Suppose that for all $n\geq 1$ and for all $x\in B_{\mu}(r)$,
\begin{equation*}
\mathcal{W}_1(\delta_x\mathcal K^n,\lambda)\leq \Lambda(x)\alpha(n),
\end{equation*}
where $\Lambda\in L^2(\lambda)$ and
\begin{equation*}
\sum_{k=1}^{\infty} \frac{\alpha(k)}{\sqrt k}<\infty.
\end{equation*}
Then condition \eqref{CLT at one point} holds with covariance matrix $4\Sigma$.
Consequently,
\begin{equation*}
\sqrt n \log_\mu\mu_n
\longrightarrow
N_m\left(0,4(D^2 F(0))^{-1}\Sigma(D^2 F(0))^{-1}\right),
\end{equation*}
in distribution, where $\Sigma:T_\mu M\to T_\mu M$ is the covariance operator determined by, for all $v\in T_\mu M$,
\begin{equation*}
\langle v,\Sigma v\rangle=\lim_{n\to\infty}\frac{1}{n}E\left[\left(
\sum_{k=1}^n\langle v,\log_\mu X_k\rangle
\right)^2
\right].
\end{equation*}
\end{theorem}
\begin{proof}
Fix non-zero $v\in T_\mu M$ and define
\begin{equation*}
g_v(x)=\frac{\Delta_U}{\norm{v}}\inpr{v}{\log_\mu x},
\end{equation*}
where
\begin{equation*}
\Delta_U=\min_{0\leq t\leq r}\operatorname{sinc}_{K_U}(t).
\end{equation*}
Then
\begin{equation*}
S_n(g_v)=\sum_{k=1}^n g_v(X_k)=\frac{\Delta_U}{\norm{v}}\sum_{k=1}^n \inpr{v}{\log_\mu X_k}.
\end{equation*}
By the Cauchy-Schwarz inequality,
\begin{equation*}
\abs{g_v(x)-g_v(y)}=\frac{\Delta_U}{\norm{v}}\abs{\inpr{v}{\log_\mu x-\log_\mu y}}\leq \Delta_U\norm{\log_{\mu}x-\log_{\mu}y}.
\end{equation*}
Let $\eta$ be the minimal geodesic from $x$ to $y$. 
Such a geodesic lies entirely in $B_{\mu}(r)$ by strong convexity, \citep{Chavel_2006}. 
Since
\begin{equation*}
\frac{d}{ds}\log_\mu(\eta(s))=D\log_\mu|_{\eta(s)}\big[\eta'(s)\big],
\end{equation*}
we obtain
\begin{equation*}
\norm{\log_{\mu}x-\log_{\mu}y}
\leq
\int_{0}^1\norm{D\log_{\mu}|_{\eta(s)}}_{\mathrm{op}}\norm{\eta'(s)}ds.
\end{equation*}
Using Lemma \ref{lem:uniform-exp-log-bounds},
\begin{equation*}
\norm{\log_{\mu}x-\log_{\mu}y}\leq \frac{1}{\Delta_U}\int_{0}^1\norm{\eta'(s)}ds
=\frac{d(x,y)}{\Delta_U}.
\end{equation*}
Therefore,
\begin{equation*}
\abs{g_v(x)-g_v(y)}\leq d(x,y),
\end{equation*}
for all $x,y\in B_{\mu}(r)$.
Since $\mu$ is the Fréchet mean, $E_\lambda[\log_\mu X_0]=0$.
Hence, $E_\lambda[g_v(X_0)]=0$.
By the $\mathcal{W}_1$-mixing assumption, Theorem 3 of \citep{jin2026centrallimittheoremsmarkov} gives
\begin{equation*}
\frac{1}{\sqrt n}\sum_{k=1}^n \inpr{v}{\log_\mu X_k} \longrightarrow N(0,\sigma_v^2),
\end{equation*}
in distribution, where
\begin{equation*}
\sigma_v^2=\lim_{n\to\infty}\frac{1}{n}E\left[\left(\sum_{k=1}^n\inpr{v}{\log_\mu X_k}\right)^2\right].
\end{equation*}

For $v,w \in T_\mu M$, define
\begin{equation*}
\inpr{v}{\Sigma w}=\lim_{n\to\infty}
\frac{1}{n}
E\left[\left(\sum_{k=1}^n \inpr{v}{\log_\mu X_k}\right)
\left(\sum_{k=1}^n \inpr{w}{\log_\mu X_k}\right)\right].
\end{equation*}
These limits define a symmetric covariance operator
$\Sigma:T_\mu M\to T_\mu M$ satisfying
\begin{equation*}
\langle v,\Sigma v\rangle=\sigma_v^2.
\end{equation*}
Hence, by the Cramér--Wold theorem,
\begin{equation*}
\frac{1}{\sqrt n}\sum_{k=1}^n \log_\mu X_k \longrightarrow N_m(0,\Sigma),
\end{equation*}
in distribution.
By Lemma \ref{lemma: bounded curvature and Hessian a.s. uniform convergence}, condition (\ref{Integrability tail-curvature condition}) holds. 
Hence, by Theorem~\ref{thm:local clt for sample Frechet mean}, the conclusion follows.
\end{proof}

\section{Application to Random Dynamical Systems}
In this section, we focus on a class of stochastic processes to which we will apply a central limit theorem.
Namely, we consider a stochastic process generated by a random dynamical system (RDS), also known as a random system of iterated functions.
This, in turn, yields a CLT for the sample Fréchet mean of a sequence generated by random dynamics.
Although random dynamical systems can be defined in the context of metric spaces, we restrict ourselves to Riemannian manifolds, since our goal is to apply the CLT from the previous section to the invariant law of an RDS.
\begin{definition}
Let $(E,\mathcal E,\eta)$ be a probability space, and let $(\xi_n)_{n\ge 1}$ be an i.i.d. sequence of $E$-valued random variables with common law $\eta$.
Let $M$ be a complete $m$-dimensional Riemannian manifold and $S\subset M$.
Let
\begin{equation*}
G:E\times S\to S,
\end{equation*}
be a measurable map. 
For a fixed $\xi\in E$, we write $G_{\xi}(x)$ to denote $G(\xi,x)$.
Define the process $(X_n)_{n\ge 0}$ by
\begin{equation*}
X_{n+1}=G_{\xi_{n+1}}(X_n),
\end{equation*}
for all $n\geq 0$, where $X_0$ is some $S$-valued random variable independent of $\{\xi_n\}_{n\geq 1}$.
We call such a process a one-step random dynamical system, or simply a random dynamical system.
\end{definition}
\begin{remark}
In the literature, random dynamical systems are often defined in a more general setting; see \citep{arnold98random}.
Specifically, a random dynamical system over a measure-preserving dynamical system $(\Omega,\mathcal{F},P,(T_t)_{t\in \mathcal{I}})$, where $\mathcal{I}$ is an abelian semigroup and $(T_t)_{t\in\mathcal{I}}$ is a family of measurable, measure-preserving maps, is a measurable map
$\Phi\colon \mathcal{I} \times \Omega \times M \to M$
that satisfies the cocycle property.
In our setting, the sequence $(\xi_n)$ is assumed to be i.i.d., which corresponds to taking $(\Omega,\mathcal{F},P)$ as the canonical product space $(E^{\mathbb{N}},\mathcal{E}^{\otimes\mathbb{N}},\eta^{\otimes \mathbb{N}})$, equipped with the left shift $T\colon E^{\mathbb{N}}\to E^{\mathbb{N}}$, defined by $T(\omega_1,\omega_2,\ldots) = (\omega_2,\omega_3,\ldots)$. 
This defines a measure-preserving dynamical system, and by defining $\Phi(0,\omega,x)=x$ and
\begin{equation*}
\Phi(n, \omega, x) = G_{\xi_n(\omega)} \circ \cdots \circ G_{\xi_1(\omega)}(x),
\end{equation*}
we recover the standard form of a random dynamical system as described in \citep{arnold98random}.
\end{remark}
The first property of such an RDS is that it is a Markov chain; see Proposition 3.1 in \citep{benaim2022markov}.
It is worth mentioning that the fact that an RDS is a Markov chain holds for Polish metric spaces; however, we restrict the model to Riemannian manifolds.
We now turn to the question of ergodicity of the model.
Specifically, when does such a process $X_n=G_{\xi_n}(X_{n-1})$ admit an invariant law?
We give a standard sufficient condition in the following lemma.
\begin{lemma}
\label{lemma: intrinsic W1 contraction random maps}
Let $M$ be a complete Riemannian manifold and $S\subset M$ be closed.
Let $(X_n)_{n\geq 0}$ be a one-step random dynamical system generated by $G_{\xi}$ such that
for some $x_0\in S$,
\begin{equation*}
E_{\eta}[d(G_\xi(x_0),x_0)]<\infty.
\end{equation*}
Define
\begin{equation*}
q=E_\eta[L(\xi)],
\end{equation*}
and assume $q<\infty$, where
\begin{equation*}
L(\xi)=\sup_{\substack{x,y\in S\\ x\neq y}}\frac{d(G_\xi x,G_\xi y)}{d(x,y)}.
\end{equation*}
Then, for all $x,y\in S$,
\begin{equation*}
\mathcal{W}_1(\delta_x\mathcal K,\delta_y\mathcal K)
\leq
q d(x,y).
\end{equation*}
Consequently, for all $\zeta,\rho\in\mathsf{Prob}_1(S)$,
\begin{equation*}
\mathcal{W}_1(\zeta\mathcal K,\rho\mathcal K)
\leq
q\,\mathcal{W}_1(\zeta,\rho).
\end{equation*}
Moreover, for every $n\geq 1$,
\begin{equation*}
\mathcal{W}_1(\zeta\mathcal K^n,\rho\mathcal K^n)
\leq
q^n \mathcal{W}_1(\zeta,\rho).
\end{equation*}
If $q<1$, then there exists a unique invariant law $\lambda$ for $\mathcal{K}$ and, therefore, $\lambda$ is ergodic.
\end{lemma}
\begin{proof}
Fix $x,y\in S$. 
For $n\geq 1$, consider the two trajectories
\begin{equation*}
X_n^x=G_{\xi_n}\circ\cdots\circ G_{\xi_1}(x),
\qquad
Y_n^y=G_{\xi_n}\circ\cdots\circ G_{\xi_1}(y),
\end{equation*}
with $X_0^x=x$ and $Y_0^y=y$.
For $n\geq 0$, let
\begin{equation*}
\mathsf C_{x,y}^{(n)}=\operatorname{Law}(X_n^x,Y_n^y).
\end{equation*}
Then $\mathsf C_{x,y}^{(n)}$ is a coupling of $\delta_x\mathcal K^n$ and $\delta_y\mathcal K^n$. Hence,
\begin{equation*}
\mathcal{W}_1(\delta_x\mathcal K^n,\delta_y\mathcal K^n)
\leq E[d(X_n^x,Y_n^y)].
\end{equation*}

On the other hand, we have
\begin{equation*}
d(X_n^x,Y_n^y) \leq \prod_{k=1}^n L(\xi_k)\,d(x,y).
\end{equation*}
Since $(\xi_k)_{k\geq 1}$ is a sequence of i.i.d. random variables, we obtain
\begin{equation*}
E[d(X_n^x,Y_n^y)] \leq E\left[\prod_{k=1}^n L(\xi_k)\right]d(x,y)
=\prod_{k=1}^nE[L(\xi_k)]d(x,y)=q^n d(x,y).
\end{equation*}
Therefore,
\begin{equation*}
\mathcal{W}_1(\delta_x\mathcal K^n,\delta_y\mathcal K^n) \leq q^n d(x,y).
\end{equation*}

Now let $\zeta,\rho\in\mathsf{Prob}_1(S)$, and let $\gamma\in\Pi(\zeta,\rho)$ be any coupling of $\zeta$ and $\rho$.
Define a probability measure $\gamma_n$ on $S\times S$ by
\begin{equation*}
\gamma_n(B)=\int_{S\times S} \mathsf C_{x,y}^{(n)}(B)\gamma(dx,dy),
\end{equation*}
for all $B\in\mathcal B(S\times S)$.
Then $\gamma_n$ is a coupling of $\zeta\mathcal K^n$ and $\rho\mathcal K^n$.
Consequently,
\begin{equation*}
\mathcal{W}_1(\zeta\mathcal K^n,\rho\mathcal K^n)
\leq \int_{S\times S} d(x',y')\,\gamma_n(dx',dy').
\end{equation*}
By the definition of $\gamma_n$,
\begin{equation*}
\int_{S\times S} d(x',y') \gamma_n(dx',dy')
= \int_{S\times S} E[d(X_n^x,Y_n^y)]\gamma(dx,dy).
\end{equation*}
Using the previous estimate,
\begin{equation*}
\int_{S\times S}E[d(X_n^x,Y_n^y)]\,\gamma(dx,dy)
\leq q^n \int_{S\times S} d(x,y)\gamma(dx,dy).
\end{equation*}
Thus,
\begin{equation*}
\mathcal{W}_1(\zeta\mathcal K^n,\rho\mathcal K^n)
\leq q^n \int_{S\times S} d(x,y)\,\gamma(dx,dy).
\end{equation*}
Taking the infimum over all couplings $\gamma\in\Pi(\zeta,\rho)$, we get
\begin{equation*}
\mathcal{W}_1(\zeta\mathcal K^n,\rho\mathcal K^n)
\leq q^n \mathcal{W}_1(\zeta,\rho),
\end{equation*}
for all $n\geq 1$.
The map defined by the action of the Markov kernel on $\mathsf{Prob}_1(S)$
$\Psi:\mathsf{Prob}_1(S)\to\mathsf{Prob}_1(S), \Psi(\zeta)=\zeta\mathcal K$, is well-defined as
\begin{align*}
\int_S d(y,x_0)\,(\zeta\mathcal K)(dy)
&=\int_S E_\eta[d(G_\xi(x),x_0)]\,\zeta(dx)
\\
&\leq \int_S E_\eta[d(G_\xi(x),G_\xi(x_0))]\,\zeta(dx)
+E_\eta[d(G_\xi(x_0),x_0)]
\\
&\leq q\int_S d(x,x_0)\,\zeta(dx)
+E_\eta[d(G_\xi(x_0),x_0)]
<\infty.
\end{align*}
Since $(S,d)$ is complete and separable, $(\mathsf{Prob}_1(S),\mathcal{W}_1)$ is complete. 
Therefore, if $q<1$, then by Banach's fixed point theorem, there exists a unique $\lambda\in\mathsf{Prob}_1(S)$ such that $\Psi(\lambda)=\lambda$, that is, $\lambda\mathcal K=\lambda$.
Finally, applying the contraction estimate with $\rho=\lambda$, and using $\lambda\mathcal K^n=\lambda$, gives
\begin{equation*}
\mathcal{W}_1(\zeta\mathcal K^n,\lambda)
= \mathcal{W}_1(\zeta\mathcal K^n,\lambda\mathcal K^n)
\leq q^n \mathcal{W}_1(\zeta,\lambda).
\end{equation*}
Since $\lambda$ is the unique invariant probability measure, it is ergodic.
This proves the claim.
\end{proof}
It is worth mentioning that the existence and uniqueness of an invariant law for an RDS may be obtained under weaker assumptions than ${E_\eta[L(\xi)]<1}$.
Specifically, it can be shown that if the map $G_{\xi}$ is a contraction on average and $E_\eta\left[\max\left\{\log d(G_{\xi}(x_0),x_0),0\right\}\right]<\infty$ for some $x_0$, then existence and uniqueness follow.
While those assumptions are also sufficient for weak convergence, they are not sufficient for $\mathcal{W}_1$-convergence; see \citep{diaconis1999iterated} or \citep{benaim2022markov}.

Let $M$ be a complete Riemannian manifold.
Given an RDS $X_{n}=G_{\xi_n}(X_{n-1})$, let $\lambda$ be an ergodic law of the induced Markov kernel, supported in a geodesic ball $B_{\mu}(r)$ with
$r<\frac{1}{2}\min\{\injr{M},\frac{\pi}{\sqrt{K_U}}\}$.
Using the uniqueness theorem of the Fréchet mean in \citep{afsari2011riemannian}, it follows that $\mu$ is the unique Fréchet mean of $\lambda$ if and only if
\begin{equation*}\label{eq: mean-zero-G} 
\int_{B_\mu(r)}E_{\eta}\left[\log_\mu\left(G_{\xi}(x)\right)\right]\lambda(dx)=0.
\end{equation*}
The same characterisation holds in a Hadamard manifold, using convexity of the distance function.

\begin{theorem}
Suppose that $M$ is a complete $m$-dimensional Riemannian manifold with
sectional curvature bounded above and below, and fix $\mu\in M$.
Let $K_U$ be an upper bound
on the sectional curvature, and suppose that
\begin{equation*}
r<\frac{1}{2}\min\left\{\injr{M},
\frac{\pi}{\sqrt{K_U}}\right\},
\end{equation*}
where $\pi/\sqrt{K_U}:=\infty$ when $K_U\leq 0$.
If $M$ is a Hadamard manifold, we take $r=\infty$, in which case
$\overline{B}_\mu(r)=M$.
Consider the RDS given by
\begin{equation*}
G_\xi:\overline{B}_{\mu}(r)\to\overline{B}_{\mu}(r)
\end{equation*}
and
\begin{equation*}
X_{n+1}=G_{\xi_{n+1}}(X_n),
\qquad n\geq 0.
\end{equation*}
Assume that, for some $x_0\in\overline{B}_{\mu}(r)$,
\begin{equation*}
E_\eta\left[d^2(G_\xi(x_0),x_0)\right]<\infty,
\end{equation*}
and that $E_\eta[L(\xi)]<1$.
Then the RDS has a unique invariant law $\lambda\in\mathsf{Prob}_1(\overline{B}_{\mu}(r))$.
Assume further that $\lambda$ has finite Fréchet function, that $X_0\sim\lambda$,
and that $\mu$ is the unique Fréchet mean of $\lambda$ and $\mu_n$ be its estimator.
Define
\begin{equation*}
F(v)=\mathcal F(\exp_\mu v),
\end{equation*}
in a neighbourhood of $0\in T_\mu M$, and assume that $D^2F(0)$ is
invertible. Then
\begin{equation*}
\sqrt n\,\log_\mu\mu_n
\longrightarrow
N_m\left(
0,
4(D^2F(0))^{-1}\Sigma(D^2F(0))^{-1}
\right),
\end{equation*}
in distribution, where $\Sigma:T_\mu M\to T_\mu M$ is the covariance
operator determined by
\begin{equation*}
\inpr{v}{\Sigma v}=\lim_{n\to\infty}\frac{1}{n}E\left[\left(\sum_{k=1}^n\inpr{v}{\log_\mu X_k}
\right)^2
\right],
\end{equation*}
for every $v\in T_\mu M$.
\end{theorem}
\begin{proof}
We show that the assumptions of Theorem \ref{thm:frechet-clt-w1-mixing} are satisfied.
By Lemma \ref{lemma: intrinsic W1 contraction random maps}, for all $\zeta,\rho\in\mathsf{Prob}_1(\overline{B}_{\mu}(r))$,
\begin{equation*}
\mathcal{W}_1(\zeta\mathcal{K}^n,\rho\mathcal{K}^n)
\leq \left(E_{\eta}[L(\xi)]\right)^n \mathcal{W}_{1}(\zeta,\rho).
\end{equation*}
Since $\lambda$ is the invariant law of $\{X_n\}_{n\geq 0}$,
\begin{equation*}
\mathcal{W}_1(\delta_x\mathcal{K}^n,\lambda\mathcal{K}^n)
\leq \left(E_\eta[L(\xi)]\right)^n \mathcal{W}_{1}(\delta_x,\lambda).
\end{equation*}
Setting $\alpha(k)=\left(E_\eta[L(\xi)]\right)^k$, we have
\begin{equation*}
\sum_{k=1}^{\infty}\frac{\alpha(k)}{\sqrt{k}}
=
\sum_{k=1}^{\infty}\frac{\left(E_\eta[L(\xi)]\right)^k}{\sqrt{k}}
\leq
\sum_{k=1}^{\infty}\left(E_\eta[L(\xi)]\right)^k.
\end{equation*}
Since $E_\eta[L(\xi)]<1$, it follows that
\begin{equation*}
\sum_{k=1}^{\infty}\frac{\alpha(k)}{\sqrt{k}}<\infty.
\end{equation*}
Next, define
\begin{equation*}
\Lambda(x)=\mathcal{W}_1(\delta_x,\lambda),
\end{equation*}
for $x\in \overline{B}_{\mu}(r)$.
Because the only coupling between $\delta_x$ and $\lambda$ is $\delta_x\otimes\lambda$, we have
$
\mathcal{W}_1(\delta_x,\lambda)=E_\lambda[d(x,Y)]$ where $Y\sim\lambda$.
By Jensen’s inequality,
\begin{equation*}
\mathcal{W}_1^2(\delta_x,\lambda)=E_\lambda[d(x,Y)]^2
\leq E_\lambda[d^2(x,Y)].
\end{equation*}
Integrating with respect to $\lambda$ gives
\begin{equation*}
\int_{M} \Lambda(x)^2\,\lambda(dx)
\leq \iint_{M\times M} d(x,y)^2\,\lambda(dx)\lambda(dy).
\end{equation*}
Moreover, the inequality
\begin{equation*}
d(x,y)^2 \leq 2\,d(x,p)^2 + 2d(y,p)^2,
\end{equation*}
for all $p,x,y\in M$, implies
\begin{align*}
\iint_{M\times M} d(x,y)^2\,\lambda(dx)\lambda(dy)
&\leq 2\int_M d(x,p)^2\,\lambda(dx)
+2\int_M d(p,y)^2\,\lambda(dy) \\
&=4\int_M d(p,z)^2\,\lambda(dz)
=4\mathcal F(p).
\end{align*}
As $\mathcal{F}(p)<\infty$ for all $p\in M$,
\begin{equation*}
\int_{M} \Lambda(x)^2\lambda(dx) < \infty.
\end{equation*}
So the assumptions of Theorem~\ref{thm:frechet-clt-w1-mixing} are satisfied and the conclusion follows.
\end{proof}

\section{Conclusions and further directions}
In this paper, we established central limit theorems for the sample Fréchet mean of stationary ergodic Markov chains taking values in manifolds.
Our first result extends the classical asymptotic framework for sample Fréchet means to the Markov-dependent setting under a central limit condition at the population Fréchet mean along with local regularity assumptions.
Under the uniqueness of the population Fréchet mean and suitable regularity conditions on the Hessian of the Fréchet function, we obtained the asymptotic normality of $\sqrt{n}\log_\mu\mu_n$.
We then provided sufficient geometric and probabilistic conditions under which the required tangent-space central limit theorem and Hessian regularity conditions hold. 
In particular, curvature bounds and support restrictions were used to control the exponential and logarithm maps, while a Wasserstein mixing-rate condition controlled the dependence structure of the Markov chain.
As an application, we considered random dynamical systems on Riemannian manifolds. 
By representing the random dynamical system as a Markov chain and imposing an average contraction condition, we verified the Wasserstein mixing-rate condition required for the central limit theorem under $\mathcal{W}_1$-mixing.

A natural direction for future work is the development of statistical inference procedures for Fréchet means under Markov dependence, including hypothesis tests and confidence regions on the manifold. 
This requires constructing consistent estimators of the limiting covariance operator, which incorporates the long-run covariance structure of the tangent-space process.
Another natural direction is to explore extensions of the central limit theorem beyond Markov dependence. 
In particular, the strong law for sample Fréchet means developed in \citep{jaffe2024fr} applies in settings more general than Markov chains.
It would therefore be interesting to investigate whether, for stationary ergodic manifold-valued sequences, such a strong law can be combined with a central limit theorem for the tangent-space process to obtain a central limit theorem for the sample Fréchet mean.

\bigskip
\noindent\textbf{Funding.}
The author is supported by a doctoral scholarship from Kuwait University.

\medskip
\noindent\textbf{Acknowledgements.}
The author would like to thank Davide Pigoli for the helpful comments and suggestions.

\bigskip
\bibliographystyle{apalike}
\bibliography{references}

\end{document}